\documentclass[12pt]{article}

\usepackage{amsmath, amsfonts, amsthm, amssymb}

\begin{document}

\title{Triangle Area Numbers and Solid Rectangular Numbers}

\author{Konstantine D. Zelator\\
Department of Mathematical Sciences\\
University of Northern Colorado\\
Greeley, CO  80639}

\maketitle

The subject matter of this paper is the area of triangles with integer side lengths and area that is a rational number.  We have the following definition:

\vspace{.25in}

\noindent {\bf Definition 1:  A triangle area number is the area number of a triangle whose sides have integer lengths and whose area is a rational number.}

As it will be established in the last section of this paper, Section 6, every triangle area number is an even integer.  Furthermore, every triangle area number is divisible by 3.  A well known subset of the set of triangle area numbers is the set of Pythagorean numbers: those are the areas of Pythagorean triangles (right triangles with integer sides).  Another interesting subset of the set of triangle area numbers is that of solid rectangular numbers (the precise defintion is given in Section 6); those numbers arise in connection with rectangular solids (or parallelepipeds) with integer edges and integer inner diagonals; more specifically, the solid rectangular numbers arise from the integer solutions to the diophantine equation (1) (see below).  The formulas obtained in {\it Result 2} of Section 5 do generate solid rectangular numbers exclusively (also see table below {\it Result 2}).  In Section 6, where a complete description of triangle area numbers is given, it is shown that the formulas of {\it Result 2} actually generate all the solid rectangular numbers and none other.  At the end of this paper (end of Section 6) we list all triangle area numbers not exceeding 999, all Pythagorean numbers and all solid rectangular numbers.  There are exactly 96 triangle area numbers not exceeding 999; thirty-one of them are Pythagorean numbers, two  are solid-rectangular numbers, but not Pythagorean, and sixty-three of them are triangle area numbers that are neither Pythagorean nor solid rectangular numbers.

\vspace{.25in}

\section{Introduction}

As it will become evident, there is a close relationship between triangles with integer side lengths and whose area $A$ is a rational number on the one hand, and the diophantine equation $x^2 + y^2 +z^2 = t^2$  (or its more general version, see (17) in Section 6) on theother hand.  Consider

\begin{equation}
x^2+y^2 + z^2 = t^2.
\label{E1}
\end{equation}

\noindent There is a specific geometric meaning behind the equation. If a quadruple of positive integers $(x_0,y_0,z_0,t_0)$ satisfies equation (\ref{E1}), then the rectangular solid whose length, width and height (twelve edges, three distinct edge length values) are the natural numbers $x_0,y_0,z_0$, will have two inner diagonals both having length equal to the natural number $t_0$.  In Section 6, we establish a process wherein solutions of (\ref{E1}) are used to generate integer-sided triangle with rational area $A$.  In order to clearly see how the relationship between such triangles and equation (\ref{E1}) arises, we must make use of the Heron-formula for the area $A$ of the triangle.  Let $a,b,c$ be positive integers, the side lengths of a triangle.

{\it Heron's formula:}  $A = \sqrt{S\cdot(s-a)(s-b)(s-c)}$ where $s = {\displaystyle \frac{N}{2}} = {\displaystyle \frac{a+b+c}{2}}$, the semi-perimeter, $N$ the perimeter of the triangle.

Squaring, we obtain $A^2 = s\cdot (s-a) \cdot (s-b) \cdot (s-c)$.  If $\theta_1,\ \theta_2, \ \theta _3$ are the degree (or radian) measures of the triangle's three angles, we know that 

$$
A = {\displaystyle \frac12} a\cdot b \cdot \sin(\theta_3) = {\displaystyle \frac12} a\cdot c\cdot \sin (\theta_2) = {\displaystyle \frac12} b\cdot c\cdot \sin (\theta_1).
$$

\noindent \framebox{\begin{tabular}{l}One can immediately observe that a triangle with integer sides $a,b,c$ will \\have rational area $A$, if and only if all three sine values $\sin(\theta_1),\ \sin(\theta_2),$\\ $ \sin(\theta_3)$ are rational numbers.\end{tabular}}

Also, any triangle with integer sides $a,b,c$ will have \underline{rational} cosine values $\cos(\theta_1),\cos(\theta_2),\cos(\theta_3)$ (regardless of whether or not $A$ is a rational number); this is an immediate consequence of the Law of Cosines

If we set

\begin{equation}
\left.\begin{array}{rcl}
-a+b+c &  = & N_1\\
\\
a-b+c & = & N_2\\
\\
a+b-c & = &  N_3 
\end{array}\right\}
\label{E2}
\end{equation}

\noindent We also obtain from solving for $a,b,c$ the following expressions:

\begin{equation}
\begin{array}{lll}
a = {\displaystyle \frac{N_2 + N_3}{2}},\ \  &\ \  b = {\displaystyle \frac{N_1 + N_3}{2}},\ \  & \ \ c = {\displaystyle \frac{N_1+N_2}{2}}
\end{array}
\label{E3}
\end{equation}

\noindent From $A^2 = s \cdot (s-a)\cdot (s-b)\cdot (s-c),\ s={\displaystyle \frac{N}{2}},\ \ N=a+b+c$, and (\ref{E3}) we obtain,

\begin{equation}
16A^2 = N \cdot N_1 \cdot N_2 \cdot N_3
\label{E4}
\end{equation}

\noindent Let $\delta$ be the greatest common divisor of the integers $a,b,c$.  And $d$ the greatest common divisor of the positive integers $N_1,N_2$, and $N_3$.  Note that $N_1,N_2$, and $N_3$ are indeed positive since $b+c > a,a+c > b, a+b > c$, these are the triangle inequalities.  In the language of number theory,

\begin{equation}
\left. \begin{array}{l}
(a,b,c) = \delta \\
\\
(N_1,N_2,N_3) = d 
\end{array}
\right\}
\label{E5}
\end{equation}

\section{The Relationship Between $\delta $ and $d$ (Irrespective of Whether the Area is Rational)}

How do integers $\delta$ and $d$ relate?  As it will be evident below, there are only two possibilities:  either $\delta =d$ or $2\delta =d$.

First note, by inspection, that (\ref{E2}) clearly shows that \underline{$\delta$ must be a divisor} \underline{of $d$}.  Also, observe that depending on the combination of parities of $a,b$, and $c$,

\vspace{.15in}

\framebox{either all three $N_1,N_2$, and $N_3$ are odd; or they are all even.}

\vspace{.15in}

In case all three $N_1,N_2$, and $N_3$ are odd, it is easy to see from (\ref{E3}) that $d$ must also be a common divisor of $a,b$, and $c$; and hence a divisor of $\delta$.  So when $N_1,N_2$, and $N_3$ are odd, and $\delta$ and $d$ must divide each other, and since they are positive, they must be equal.  However, when $N_1,N_2$, and $N_3$ are all even, the situation becomes significantly more complicated.  Without elaborating on the details, an analysis modulo 4 demonstrates that the following possibilities or combinations occur:

\noindent (Parenthetically we note here that in the above observations, an underlying assumption was made- a fact very well known in number theory: every common divisor of any number of integers, must divide their greatest common divisor.)

Below, we include all possibilities, including the case where $N_1,N_2$, and $N_3$ are odd.  This analysis is the result of equations (\ref{E2}) and (\ref{E3}), and considerations modulo 4.  In total, we have seven broad occurrences or combinations:

\vspace{.25in}

\noindent \framebox{\begin{tabular}{l}{\bf \underline{Case 1:}}  $a,b,c$ are odd:  $N_1,N_2$, and $N_3$ are odd, $\delta = d$.\\
{\bf \underline{Case 2:}}  One of $a,b,c$ is even and the other two odd:  one of $N_1,N_2$, and \\$N_3$ is congruent to 0 modulo 4; the other two are congruent to 2(mod 4).\\  In this case, $d = 2\delta$.\\
{\bf \underline{Case 3:}} Two of $a,b,c$ are even, the other odd:  $N_1,N_2,N_3$ are odd,\\ and $\delta = d$.\\
{\bf \underline{Case 4:}} All three of $a,b,c$ are even, with either all three being congruent to\\
2(mod 4); or  with two of them congruent to 0(mod 4), the third one 2(mod 4).\\  In either 
(sub)-case, all three $N_1,N_2,N_3$ must be congruent to 2(mod 4); $\delta =d$.\\
{\bf \underline{Case 5:}} Two of $a,b,c$ are congruent to 2(mod 4), the third one congruent to\\ 0(mod 4): Then, all three $N_1,N_2,N_3$ are multiples of 4; and $d = 2 \delta$.\\
{\bf \underline{Case 6:}} All three $a,b,c$, are multiples of 4; and so must $N_1,N_2$, and $N_3$ be.\\
And with the exponents $\alpha_{i_1},\alpha_{i_2}, \alpha_{i_3}$ satisfying $\alpha_{i_1} \geq \alpha_{i_2} \geq \alpha_{i_3}$, with at least one \\
inequality sign being strict (so that not all $\alpha_{i_1},\alpha_{i_2},\alpha_{i_3}$ are equal), where\\ $(i_1,i_2,i_3)$ is a permutation of $(1,2,3)$ and $2^{\alpha_{i_1}},2^{\alpha_{i_2}},2^{\alpha_{i_3}}$, are the highest powers\\ of 2 dividing $N_{i_1},N_{i_2},N_{i_3}$, respectively.  In this case $d = 2 \delta$.\\
{\bf \underline{Case 7:}} All three $a,b,c$, are divisible by 4; and thus, so are
 $N_1,N_2$, and $N_3$.\\  And with $N_1,N_2,N_3$ being divisible by the same highest power
 of 2 (that is,\\ $\alpha_1 = \alpha_2 = \alpha_3$ in the notation of Case 6). In this case $\delta = d$.
\end{tabular}}

\vspace{.25in}

\begin{center}
{\bf \underline{Numerical Examples}}
\end{center}

\vspace{.15in}

\begin{enumerate}
\item[1)] {\bf \underline{For Case 1:}}  Any three odd natural numbers $a,b,c$, will do, as long as $a+b > c$, $a+c>b$, and $b+c >a$.
\item[2)] {\bf \underline{For Case 2:}}  $a = 6,b=5,c=9,N_1=8,N_2=10,N_3=2;\ \delta =1,\d=2$
\item[3)] {\bf \underline{For Case 3:}} $a=6,b=12,c=9,N_1=15,N_2=3,N_3=9;\ \delta=d=3$
\item[4)] {\bf \underline{For Case 4:}} $a=6,b=14,c=18,N_1=26,N_2=10,N_3 =2;\ \delta =d=2$
\item[5)] {\bf \underline{For Case 5:}} $a=34,b=38,c=44,N_1 = 48,N_2 =40,N_3=28;\ \delta =4,d=8$
\item[6)] {\bf \underline{For Case 6:}} $a=34, b=38, c=44, N_1 = 48,N_2=40,N_3=28;\ \delta=4,d=8,\ \alpha_1 =3,\alpha_2=5,\alpha_3=4$
\item[7)] {\bf \underline{For Case 7:}} $a = 24,b=28,c=36,N_1=40,N_2=32,N_3=16; \ \delta =d=8,\ \alpha_1=\alpha_2=\alpha_3$
\end{enumerate}

\section{The Primitive Case and Result 1}

According to (\ref{E5}), we have $\left( {\displaystyle \frac{a}{\delta},\frac{b}{\delta},\frac{c}{\delta}}\right) = 1 = \left( {\displaystyle \frac{N_1}{d},\frac{N_2}{d},\frac{N_3}{d}} \right)$.  Furthermore, note that since the perimeter $N$ is equal to $N=a+b+c$, it is also true that $N=N_1+N_2+N_3$ (by adding the three equations in (\ref{E2}); or in (\ref{E3})).  In what is to follow, we will make the assumption below:

\begin{equation}
\left( {\displaystyle \frac{N_1}{d}\cdot \frac{N_2}{d} \cdot \frac{N_3}{d} ,\frac{N}{d}}\right)= 1 
\label{E6}
\end{equation}

\noindent In other words, we will assume that $\left({\displaystyle \frac{N}{d}}\right)$ is relatively prime to each of the integers ${\displaystyle \frac{N_1}{d},\frac{N_2}{d},\frac{N_3}{d}}$.

\vspace{.15in}

\noindent\framebox{\begin{tabular}{l}{\bf {\it \underline{Remark 1:}}}  {\it It can be proven that} (\ref{E6}) {\it implies} $\left({\displaystyle \frac{a}{\delta}\cdot \frac{b}{\delta}\cdot \frac{c}{\delta}, \frac{N}{d}} \right) = 1$ {\it in \underline{Cases}}\\
{\it 1, 2, 3, 4, 5, and 7, listed in the previous section}; {\it also in \underline{Case 6}, except}\\{\it when} $\alpha_{i_1}=\alpha_{i_2} > \alpha_{i_3}$ {\it or when}
$\alpha_{i_1}> \alpha_{i_2} = \alpha_{i_3}$;  {\it in which sub-cases}\\ {\it (of \underline{Case 6})}, {\it we have} $\left({\displaystyle \frac{a}{\delta} \cdot \frac{b}{\delta} \cdot \frac{c}{\delta},\frac{N}{\delta}}\right) = 2$. 
{\it We leave this as an exercise}\\ {\it for the reader to prove}.
\end{tabular}}

\vspace{.25in}

Now, let us revisit equation (\ref{E4}).  We have,

\begin{equation}
16A^2 = (4A)^2 = d^4 \cdot \left({\displaystyle \frac{N}{d}}\right) \cdot \left({\displaystyle \frac{N_1}{d}}\right) \cdot \left( {\displaystyle \frac{N_2}{d}} \right)\cdot \left( {\displaystyle \frac{N_3}{d}}\right) 
\label{E7}
\end{equation}

\noindent The right-hand side of (\ref{E7}) is an integer, while the left-side is the square of a rational number.  It is a standard exercise in elementary number theory to prove that if the square of a rational number is an integer, then the rational number itself must be an integer.  In our situation this says that if $4A$ is a rational number (which is equivalent to saying that $A$ is rational), then in fact, $4A$ must be an integer.  Hence, the left side of (\ref{E7}), will be an integer square when $A$ is rational.  Now apply condition (\ref{E6}):  Since the two factors ${\displaystyle \frac{N}{d}}$ and ${\displaystyle \frac{N_1}{d},\frac{N_2}{d},\frac{N_3}{d}}$ are relatively prime and the product is an integral square, each of them must be an integer square.  (Another standard exercise in elementary number theory.)  Thus, we must have

\begin{equation}
{\displaystyle \frac{N_1}{d} \cdot \frac{N_2}{d} \cdot \frac{N_3}{d}} = k^2
\label{E8}
\end{equation}

\noindent and 

$$
{\displaystyle \frac{N}{d}} = t^2
$$

\noindent for some integers $k$ and $t$.  But we know that $N_1 +N_2 + N_3 = N$ (refer to (\ref{E2}) or (\ref{E3})), so that combined with ${\displaystyle \frac{N}{d}}= t^2$ we obtain,

\begin{equation}
{\displaystyle \frac{N_1}{d} + \frac{N_2}{d} + \frac{N_3}{d}} = t^2
\label{E9}
\end{equation}

\noindent Let $d_{12}$ be the greatest common divisor of ${\displaystyle \frac{N_1}{d}}$ and ${\displaystyle \frac{N_2}{d}}$; likewise let $d_{13}$ and $d_{23}$ be the greatest common divisors of ${\displaystyle \frac{N_1}{d}}$ and ${\displaystyle \frac{N_3}{d}}$, and ${\displaystyle \frac{N_2}{d}}$, and ${\displaystyle \frac{N_3}{d}}$; respectively. A cursory observation reveals that these three natural numbers must be mutually relatively prime (that is any two of them must be relatively prime).  Why? Because if a prime divisor divided any two of them, it would have to divide all three ${\displaystyle \frac{N_1}{d},\frac{N_2}{d}}$, and ${\displaystyle \frac{N_3}{d}}$ thus violating the condition $\left( {\displaystyle \frac{N_1}{d}, \frac{N_2}{d}, \frac{N_3}{d}} \right) = 1$.

We can now write,

\begin{equation}
\left.
\begin{array}{l}
\left( d_{12},d_{13} \right) = \left( d_{12},d_{23}\right) = \left(d_{13}, d_{23}\right) =1,\ {\displaystyle \frac{N_1}{d}} = d_{12} \cdot d_{13} \cdot M_1\\
\\
{\displaystyle \frac{N_2}{d}} = d_{12} \cdot d_{23} \cdot M_2, {\displaystyle \frac{N_3}{d}} = d_{13} \cdot d_{23} \cdot M_3,\\
\\
{\rm for \ some\ positive\ integers}\\
\\
 M_1,M_2,M_3\ {\rm with}\ \left(M_1, M_2\right)  = \left(M_1,M_3\right)= \left(M_2,M_3\right) = 1
\end{array} \right\}
\label{E10}
\end{equation}

If we combine the formulas of (\ref{E10}) with (\ref{E8}) we obtain

$$
d^2_{12} \cdot d^2_{13} \cdot d^2_{23} \cdot M_1 \cdot M_2 \cdot M_3 = k^2
$$

\noindent The latter equation clearly shows that the product $M_1 \cdot M_2 \cdot M_3$ must be an integer square 

\begin{equation}
M_1 = k^2_1,\ \ \ \ M_2 = k^2_2,\ \ \ \ M_3 = k^2_3 
\label{E11}
\end{equation}

\noindent for natural numbers $k_1,k_2$, and $k_3$.

Combining (\ref{E11}) with the formulas in the second line of (\ref{E10}), and with (\ref{E9}) we arrive at the following necessary condition:

\vspace{.15in}

\hspace{.5in}\framebox{$d_{12}\cdot d_{13} \cdot k^2_1 + d_{12} \cdot d_{23} \cdot k^2_2 + d_{13} \cdot d_{23} \cdot k^2_3 = t^2$}\hfill (12)

\setcounter{equation}{12}

\vspace{.15in}

\noindent An immediate consequence of this necessary condition is the following result:  
\newpage
{\bf As part of the hypothesis of Result 1, assume condition (\ref{E6})}

\vspace{.15in}

\noindent\framebox{$\!\!\!\!\!$ 
\begin{tabular}{l} {\bf Result 1:} {\it Let} $a,b,c$ {\it be the integer side-lengths of a triangle; and},\\ 
$N_1 = -a + b + c,N_2 = a - b + c,\ N_3 = a + b - c,\ d = (N_1,N_2,N_3),$\\
$d_{12} = \left( {\displaystyle \frac{N_1}{d},\frac{N_2}{d}} \right),
d_{13} = \left( {\displaystyle \frac{N_1}{d},\frac{N_3}{d}}\right),d_{23} = \left( {\displaystyle \frac{N_2}{d},\frac{N_3}{d}} \right)$.  {\it Also let} $i_1,i_2,i_3$\\ {\it be a permutation} {\it of} $(1,2,3)$.\\
\noindent \begin{tabular}{ll} (a) & {\it There is no triangle with all three} $a,b,c,$ {\it being odd or with },\\
& {\it two of them even the other odd; and with the area} $A$\\
&  {\it being a rational number.  Moreover, there exists no triangle with }\\
& {\it the area being a rational number  and the three integers} $N_1,N_2,N_3$\\
& {\it being divisible by the same highest power of} $2$.\\
(b) &  {\it There exists no triangle with rational area}, ${\displaystyle \frac{N_{i_1}}{d}}$ {\it even},\\
& ${\displaystyle \frac{N_{i_2}}{d}}$ {\it and} ${\displaystyle \frac{N_{i_3}}{d}}$ {\it both odd}, {\it and} \\
&$d_{i_1i_2} \equiv 1 \equiv d_{i_1i_3}$ (mod 4); {\it or with} $d_{i_1i_2}  \equiv 3 \equiv d_{i_1i_3}$ (mod 4).\end{tabular}
\end{tabular}}

\vspace{.15in}

\noindent \framebox{\begin{tabular}{ll} (c) & {\it If}  $d_{12} = d_{23} = d_{13} =1$, {\it then there is no triangle with rational area} $A$.
\end{tabular}}

\vspace{.15in}

\noindent{\underline{\bf Note:}} In part (b), since ${\displaystyle \frac{N_{i_2}}{d},\frac{N_{i_3}}{d}}$, are assumed to be odd, it follows from (\ref{E10}) that all three $d_{i_1i_2},d_{i_1i_3} , d_{i_2i_3}$ must be odd.

\vspace{.15in}

\noindent {\underline{\bf Proof:}}

\begin{enumerate}
\item[(a)] Since $a,b,c$ are odd or two of them even the other odd; or since $N_1,N_2,N_3$ are divisible by the same highest power of $2$; the integers ${\displaystyle \frac{N_1}{d},\frac{N_2}{d},\frac{N_3}{d}}$ are odd, we argue by contradiction:  If such a triangle existed, then the necessary condition (12) would be satisfied.  Because ${\displaystyle \frac{N_1}{d},\frac{N_2}{d},\frac{N_3}{d}}$ are odd, (\ref{E10}) and (\ref{E11}) imply that $d_{12},d_{13}, d_{23},k_1,k_2$ and $k_3$ must all be odd; thus (12) implies that $t$ must also be odd, since $t^2$ is the sum of three odd integers.  But the square of any odd integer must be congruent to 1 modulo 4.  Hence (12) would imply, $d_{12} \cdot d_{13} + d_{12}\cdot d_{23}+d_{13}\cdot d_{23} \equiv 1$ (mod 4).  This congruence though is impossible, no matter what the three odd numbers $d_{12},d_{13},d_{23}$ are.  Why?  Because $d_{12}\cdot d_{13} + d_{12} \cdot d_{23} + d_{13} \cdot d_{23}$ is always congruent to 3 modulo 4.  We leave it to the reader to verify that this is the case in each of the following combinations:

\noindent $d_{12}\equiv d_{13} \equiv d_{23}$ (mod 4) (so either all three are $\equiv 1$ (mod 4); or $\equiv 3$ (mod 4); two of $d_{12},d_{13},d_{23}$ are congruent to 1 (mod 4); the other one congruent 3 (mod 4).  Or vise versa.
\item[(b)]  ${\displaystyle \frac{N_{i_2}}{d}, \frac{N_{i_3}}{d}}$ are odd,  it follows from  (\ref{E10}) that $d_{i_1i_2},d_{i_1i_3}, d_{i_2i_3}$, are all odd; and from (\ref{E10}) and  (\ref{E11})  that $k_{i_2},k_{i_3}$ are both odd while $k_{i_1}$ is even by (\ref{E10}), (\ref{E11}), and the assumption that ${\displaystyle \frac{N_{i_1}}{d}}$ is even.  Again, the argument is by contradiction: since $k^2_{i_1} \equiv 0$ (mod 4); and since the integer $t$ in (12) must be even, (12) would imply, $d_{i_2i_3} \cdot \left(d_{i_1i_2} + d_{i_1i_3}\right) \equiv 0$ (mod 4); and since $d_{i_2i_3}$ is odd, the last congruence would imply $d_{i_1i_2} + d_{i_1i_3} \equiv 0$ (mod 4), which is impossible, since the hypothesis in part (b) implies $d_{i_1i_2} + d_{i_1i_3} \equiv 2$ (mod 4).
\item[(c)] Since $d_{12} = d_{23} = d_{13} = 1$, (12) implies $k^2_1 + k^2_2 + k^2_3 = t^2$.  But the integers $k_1,k_2,k_3$ are mutually relatively prime; hence either all three are odd or one of them is even, the other two odd.  Thus $k^2_1 + k^2_2 +k^2_3 \equiv 3$ or 2(mod 4); but $t^2 \equiv 1$ or 0(mod 4), and so we have a contradiction.
\end{enumerate}

\vspace{.15in}

\noindent \framebox{\begin{tabular}{l} {\underline{\bf An Immediate Consequence of Result 1(a):}} If an integer-sided\\
triangle satisfying (\ref{E6}) has rational area $A$, then $A$ must be, in fact, an\\ 
\underline{even integer}. To see why, consider the numbers $N_1,N_2,N_3$, as defined in (\ref{E2}).\\  
Clearly either all three are odd or all three are even.  If they are all odd,\\ 
then that means, either the three sides $a,b,c$ are odd; or two are even, the\\
third odd; but then, according to Result 1(a), such a triangle cannot have\\
rational area. Therefore, if $A$ is rational, all three $N_1,N_2,N_3$ must be\\
even, and consequently so must be the perimeter $N=N_1 + N_2 + N_3$.\\
According to (\ref{E4}), $A = \sqrt{\frac{N \cdot N_1\cdot N_2 \cdot N_3}{16}}$;\\
however, the product $N\cdot N_1 \cdot N_2 \cdot N_3$ will be divisible by 16;\\
thus, $A$ will be the square root of the integer; so if $A$ is rational,\\ 
it must be an integer. To see that $A$ must be even, observed that $A$ would\\
be odd only when each of $N_1,N_2,N_3,N$ is exactly divisible by $2$; so that\\
the product $N_1N_2N_3N$ is exactly divisible by 16;  but that would mean, in\\
particular, that the highest power of $2$ dividing all three $N_1,N_2,N_3$ would\\
be $2^1$; this is precluded by Result 1a though.\end{tabular}}

\section{The General Solution to Diophantine Equation (1):}  This can be found in some number theory books.  The derivation of the general solution (\ref{E1}) can be found in W. Sierpinski's book \underline{Elementary Theory of} \underline{Numbers} see [1].  \underline{Every} solution $(x,y,z,t)$, in positive integers $x,y,z$, and $t$, of equation (\ref{E1}) must be of the form,

\vspace{.15in}

\noindent\framebox{\begin{tabular}{l} $x = {\displaystyle \frac{l^2 + m^2 -n^2}{n}}, y = 2l,\ z = 2m,\ t = {\displaystyle \frac{l^2+m^2+n^2}{n}}$ where $m,n,l$\\ are positive integers, 
$n$ is a divisor of $l^2 + m^2$ and $n < \sqrt{l^2+m^2}$
\end{tabular} \hspace{.1in} (13)}
\setcounter{equation}{13}
\vspace{.15in}  

\noindent And conversely, any quadruple $(x,y,z,t)$ that satisfies (13) must be a solution in positive integers $x,y,z,t$, of equation (\ref{E1}).  Below, we make the following remarks or observations:

\vspace{.15in}

\noindent \underline{\bf Observation 1:}  If $(x,y,z,t)$ is a (positive) integer solution to (\ref{E1}), note that an argument modulo 4 shows that at \underline{least two} of $x,y,z,t$ must be even integers.

\noindent \underline{\bf Observation 2:}  If $(x,y,z,t)$ is a solution to (\ref{E1}) with $y$ and $z$ even, then the numbers $l = {\displaystyle \frac{y}{2}},\ m={\displaystyle \frac{z}{2}}, n= {\displaystyle \frac{t-x}{2}}$ are clearly uniquely determined; which shows that every solution of the equation (\ref{E1}), in positive integer $x,y,z,t$ with $y,z$ even, is obtained exactly once by the use of formulas (13).

\noindent \underline{\bf Observation 3:}  In order to eliminate solutions of (\ref{E1}) with interchanged unknowns, we may reject pairs, $l,m$ for which $m > l$ and take only these $n$ for which the numbers $x$ are odd.  Thus, eliminating also all the solutions of (\ref{E1}) for which $x,y,z$, and $t$ are even.  To include them again, it is sufficient to multiply each of the solutions with odd $x$ by the powers of 2 successively.

Below, we present a sample of fifteen solutions to (\ref{E1}), with $m \leq l$ and $x$ odd.

\newpage

\begin{center}
{\bf A Sample of Solutions (with $m \leq l$ and $x$ odd)}
\end{center}
$$
\begin{array}{|c|c|c|c|c|c|c|c|}
\hline
\ \ \ L\ \ \  &\ \ \  m\ \ \  &\ \ \  l^2 +m^2 \ \ \  &\ \ \ n\ \ \ & \ \ 
 \  x\ \ \  & \ \ \ y\ \  \ &\ \  \ Z \ \ \ & \ \ \ t\ \ \ \\ \hline
1 & 1 & 2 & 1 & 1 & 2 & 2 & 3 \\ \hline
3 & 3 & 18 & 3 & 3 & 6 & 6 & 9\\ \hline
5 & 1 & 26 & 1 & 25 & 10 & 2 & 27 \\ \hline
5 & 3 & 34 & 1 & 33 & 10 & 6 & 35 \\ \hline
5 & 5 & 50 & 1 & 49 & 10 & 10 & 51\\ \hline
6 & 4 & 52 & 1 & 51 & 12 & 8 & 53 \\ \hline
7 & 3 & 58 & 1 & 57 & 14 & 6 & 59 \\ \hline
7 & 7 & 98 & 1 & 97 & 14 & 14 & 99 \\ \hline
8 & 2 & 68 & 1 & 67 & 16 & 4 & 69 \\ \hline
9 & 3 & 90 & 5 & 13 & 18 & 6 & 23\\ \hline
9 & 3 & 90 & 6 & 9 & 18 & 6 & 21 \\ \hline
9 & 5 & 106 & 2 & 51 & 18 & 10 & 55 \\ \hline
9 & 7 & 130 & 10 & 3 & 18 & 14 & 23 \\ \hline
9 & 9 & 162 & 1 & 161 & 18 & 18 & 163\\ \hline
9 & 9 & 162 & 9 & 9 & 18 & 18 & 27 \\ \hline
\end{array}
$$

\section{Result 2 and Generating Integer-sided Triangles with Rational Area:}

\vspace{.15in}

We state the result:

\vspace{.15in}

\framebox{\begin{tabular}{l} \underline{\bf Result 2:} {\it If} $(x,y,z,t)$ {\it is any solution, in positive integers} $x,y,z$,\\
 {\it and} $t$, {\it to equation} (\ref{E1}), {\it then any triangle with integer sides} $a,b,c$,\\ 
{\it defined below, has \underline{rational area}};\\
\\
$a={\displaystyle \frac{D(y^2+z^2)}{2}}, b = {\displaystyle \frac{D(x^2+z^2)}{2}},
 c = {\displaystyle \frac{D(x^2+y^2)}{2}},$\\
\\
{\it where} $D$ {\it is a positive integer, if all three} $x,y,z$ {\it are even;  and} $D$\\
 {\it is an even positive} {\it integer, if two of} $x,y,z$ {\it are even, the third one odd.}\\
\underline{\it Note:} {\it Recall that any solution to (\ref{E1}), requires that at} \\ 
{\it \underline{least two} of} $x,y,z$ {\it must be even. Also note that when one of} \\
$x,y,z$ {\it is odd,} $D$ {\it must be even, otherwise two of} $a,b,c$,  \\
{\it would not be integers}.
\end{tabular}}

\vspace{.15in}

\noindent \underline{\bf Proof:}  We compute four quantities in terms of $x,y,z$, and $t$:

$s = {\displaystyle \frac{a+b+c}{2}} = {\displaystyle \frac{D \cdot (x^2+y^2+z^2)}{2}} = {\displaystyle \frac{D \cdot t^2}{2}}$, since $x^2 + y^2 + z^2 = t^2$ according to hypothesis.

Next, we compute $s - a = {\displaystyle \frac{D\cdot(x^2+y^2+z^2)}{2}} = {\displaystyle \frac{D \cdot (y^2+z^2)}{2}} = {\displaystyle \frac{D\cdot x^2}{2}}$.

Likewise $s - b = {\displaystyle \frac{D\cdot (x^2 +y^2 + z^2)}{2}} - {\displaystyle \frac{D\cdot (x^2+z^2)}{2}} = {\displaystyle \frac{D \cdot y^2}{2}}$.

And $ s- c = {\displaystyle \frac{D\cdot (x^2 + y^2 + z^2}{2}} - {\displaystyle \frac{D\cdot (x^2+y^2)}{2}} = {\displaystyle \frac{D\cdot z^2}{2}}$.

Now apply Heron's formula:  $A = \sqrt{s\cdot (s-a)\cdot (s-b) \cdot (s-c)};\  A = {\sqrt{\frac{D^4 \cdot x^2 \cdot y^2 \cdot z^2 \cdot t^2}{2^4}}} = {\displaystyle \frac{D^2 \cdot x \cdot y \cdot z \cdot t}{4}}$, which proves that $A$ is a rational number, in fact an integer since at least two of the $x,y,z$ are even.

\vspace{.15in}

\noindent \framebox{ \begin{tabular}{l}
{\it \underline{\bf Remark 2:} Note that if in (12), each of the numbers} $d_{12} \cdot d_{13}, d_{12}$\\
$\cdot d_{23},d_{13}\cdot d_{23}$ {\it is an integer square, then by starting with an integer}\\
{\it -sided triangle with rational area, one can produce a solution to (\ref{E1})}.\\
 {\it Also, since} $d_{12}, d_{13}, d_{23}$, {\it are pairwise relatively}
 {\it prime, then all three}\\ 
{\it products} $d_{12} \cdot d_{13}, d_{12}\cdot d_{23}, d_{13} \cdot d_{23}$ {\it will be}
{\it perfect squares, if and only if}\\
{\it each of the three numbers} $d_{12}, d_{13},d_{23}$ {\it is a square itself; the smallest }\\ 
{\it such values are (without loss of generality)} $d_{12}=4$ {\it and} $d_{13} = d_{23} =1$.\\
{\it This can happen with} $k_3$ {\it odd in (12); more specifically, if we}\\
 $k_1 = k_2 = k_3 =1$ {\it and} $t = 3$ {\it then (12) is satisfied.  A calculation shows}\\
{\it  that (see (\ref{E11}) and (\ref{E10})), that if we also take} $d=2$, {\it then}\\ 
$N_1=8, N_2 =8, N_3 = 2$; {\it so that  (from (\ref{E3}))} $a = 8, b = 5, c=5$, {\it and the}\\ {\it area is} $A = 12$.
\end{tabular}}

\vspace{.15in}

\noindent Below, we use the table of sample solutions of equation (\ref{E1}) (found in the previous section) to  generate integer-sided triangles with rational area.  For this we use the formulas, found in Result 2, that generate the three sides $a,b,$ and $c$.  We take \framebox{$D = 2$},  since in the table of sample solution to (\ref{E1}), $x$ is odd.  As it will become clear, in the next section (Section 6), all the numbers listed below are {\bf solid rectangular numbers}.

\vspace{.15in}

\noindent$\begin{array}{|c|c|c|c|c|c|c|l|}
\hline
x & y & z & t & a & b & c & \begin{array}{rcl}A{\rm (area)}&  = & {\displaystyle \frac{D^2 \cdot x \cdot y \cdot z \cdot t}{4}}\\ & = &x \cdot y \cdot z \cdot t (D = 2)\end{array} \\ \hline
1 & 2 & 2&  3 & 8 & 5 & 5 & 12 = 2^2 \cdot 3\\ \hline
3 & 6 & 6 & 9 & 72 & 45 & 45 & 972=2^2\cdot 3^5\\ \hline
25 & 10 & 2 & 27 & 104 & 627 & 725 & 13,500 = 2^2 \cdot 3^3 \cdot 5 ^3 \\ \hline
33 & 10 & 6 & 35 & 136 & 1,125 & 1,189 & 69,300 = 2^2 \cdot 3^2 \cdot 5 \cdot 7 \cdot 11\\ \hline
49 & 10 & 10 & 51 & 200 & 2,501 & 2,501 & 249,900 = 2^2 \cdot d \cdot 5^2 \cdot 7^2 \cdot 17 \\ \hline
51 & 12 & 8 & 53 & 208 & 2,665 & 2,745 & 259,488 = 2^5 \cdot 3^2 \cdot 17 \cdot 53 \\ \hline
57 & 14 & 6 & 59 & 232 & 3,285 & 3,445 & 282,492 = 2^2 \cdot 3^2 \cdot 7 \cdot 19 \cdot 59 \\ \hline 
97 & 14 & 14 & 99 & 392 & 9,605 & 9,605 & 1,882,188 = 2^3 \cdot 3^4 \cdot 7^2 \cdot 11 ^2 \\ \hline 
67 & 16 & 4 & 69 & 272 & 4,505 & 4,745 & 295,872 = 2^6 \cdot 3 \cdot 23 \cdot 67\\ \hline
13 & 18 & 6 & 23 & 360 & 204 & 493 & 32,292 = 2^2 \cdot 3^3 \cdot 13 \cdot 23 \\ \hline
9 & 18 & 6 & 21 & 360 & 81 & 405 & 20,412 = 2^2 \cdot 3^6 \cdot 7\\ \hline
51 & 18 & 10 & 55 & 424 & 2,701 & 2,925 & 504,900 = 2^2 \cdot 3^3 \cdot 5^2 \cdot 11 \cdot 17 \\
\hline
3 & 18 & 14 & 23 & 520 & 205 & 333 & 17,388 = 2^2\cdot 3^3 \cdot 7 \cdot 23\\ \hline
161 & 18 & 18 & 163 & 648 & 26,245 & 26,245 & 8,502,732 = 2^2 \cdot 3^4 \cdot 7 \cdot 23 \cdot 163\\ \hline
9 & 18 & 18 & 27 & 648 & 415 & 415 & 78,732 = 2^2 \cdot 3^9\\ 
\hline
\end{array}
$

\vspace{.25in}

\section{ A complete description of triangle area numbers and solid rectangular numbers.}  

In this section we drop condition (6) of {\bf Section 3}.  The purpose of (6) was to simplify the situation in order to establish a connection with equation (\ref{E1}).  So, we will only assume the validity of equations (\ref{E1}) through (\ref{E5}) which were established in {\bf Section 1}.

Let us reconsider the integers ${\displaystyle \frac{N_1}{d}, \frac{N_2}{d}, \frac{N_3}{d}}$, and ${\displaystyle \frac{N}{d}}$.  Because of (\ref{E5}), if $D_1$ is the greatest common divisor of ${\displaystyle \frac{N_1}{d}}$ and ${\displaystyle \frac{N}{d}}$, then $D_1$ must be relatively prime to both ${\displaystyle \frac{N_2}{d}}$, and ${\displaystyle \frac{N_3}{d}}$; for if, say $D_1$ had a prime divisor in common with ${\displaystyle \frac{N_2}{d}}$, then (recall that) ${\displaystyle \frac{N_1}{d} + \frac{N_2}{d} + \frac{N_3}{d} = \frac{N}{d}}$ would imply that that prime divisor would also be a divisor of ${\displaystyle \frac{N_3}{d}}$, contrary to the fact that $\left( {\displaystyle \frac{N_1}{d},\frac{N_2}{d},\frac{N_3}{d}} \right) = 1$.  Likewise, by applying the same reasoning, we see that $D_2$, the greatest common divisor of ${\displaystyle \frac{N_2}{d}}$ and ${\displaystyle \frac{N}{d}}$, must be relatively prime to both ${\displaystyle \frac{N_1}{d}}$ and ${\displaystyle \frac{N_3}{d}}$.  And if $D_3 = \left( {\displaystyle \frac{N_3}{d}, \frac{N}{d}}\right)$, we must have $\left( D_3,{\displaystyle \frac{N_1}{d}}\right) = 1 = \left( D_3,{\displaystyle \frac{N_2}{d}}\right)$.  Next, let $d_{12} = \left( {\displaystyle \frac{N_1}{d}, \frac{N_2}{d}} \right),\ d_{13} = \left( {\displaystyle \frac{N_1}{d},\frac{N_3}{d}}\right), \ d_{23} = \left( {\displaystyle \frac{N_2}{d}, \frac{N_3}{d}} \right)$.  Again as before, by virtue of ${\displaystyle \frac{N}{d} = \frac{N_1}{d} + \frac{N_2}{d} + \frac{N_3}{d}}$ and $\left( {\displaystyle \frac{N_1}{d},\frac{N_2}{d},\frac{N_3}{d}} \right) = 1$, we conclude that the three integers $d_{12},d_{13},d_{23}$, are mutually relatively prime.  We can write,

\vspace{.25in}

\noindent \framebox{\begin{tabular}{l} ${\displaystyle \frac{N_1}{d}} = D_1d_{12}d_{13} M_1, {\displaystyle \frac{N_2}{d}} = D_2 d_{12}d_{23} M_2,\ {\displaystyle \frac{N_3}{d}} = D_3 d_{13}d_{23}M_3$, and\\
\\
 ${\displaystyle \frac{N}{d}} = D_1D_2D_3 M$;
where the positive integers\\ 
\\
$D_1,\ D_2,\ D_3,\ d_{12},\ d_{13},\ d_{23},\ M_1,\ M_2,\ M_3$, and $M$ satisfy the \\
\\
following seven {\bf coprimeness  conditions}:\\
\\
\begin{tabular}{rlrrr}
(i) & $(D_1,\ D_2\ d_{23}\ M_2) = (D_1,\ D_3\ d_{23}\ M_3) = 1$ &\ \ \ \ \ \  & \\
(ii) & $(D_2,\ D_1\ d_{13} \ M_1) = (D_2,\ D_3 \ d_{13}\ M_3) = 1$ & & \\
(iii) & $(D_3,\ D_1\ d_{12}\ M_1) = (D_3,\ D_2\ d_{12} M_2) = 1$ & &  \\
(iv) & $(d_{12},\ d_{13}) = (d_{12},\ d_{23}) = (d_{13},\ d_{23})=1$ & &\hspace{.85in} (14)\\
(v) & $(M,\ d_{12} \ d_{13}\ d_{23}\ M_1 \ M_2 \ M_3) = 1$ & & \\
(vi) & $(D_1\ D_2\ D_3,\ d_{12}\ d_{13},\ d_{23}) = 1$ & &  \\
(vii) & $(M_1,\ M_2) = (M_1,\ M_3) = (M_2,\ M_3)=1$ & & 
\end{tabular}\end{tabular}}

\vspace{.25in}

\noindent According to (\ref{E4}), $16A^2 = N \cdot N_1\cdot N_2\cdot N_3$; using the four formulas in (14) we obtain,
\setcounter{equation}{14}

\begin{equation}
16A^2 = d^4\cdot D^2_1D^2_2D^2_3d^2_{12}d^2_{13}d^2_{23} \cdot M_1 M_2 M_3 M \label{E15}
\end{equation}

\noindent Since $4A$ is a rational number and the right-hand side of (\ref{E15}) is an integer, $(4A)^2$ must be an integer; which means that $4A$ must actually be an integer (we have already come across this argument in Section 3 - see below equation (\ref{E7})).  Therefore $(4A)^2$ is an integer square so that (\ref{E15}) implies that $M_1M_2M_3M$ must be an integer square as well:  but according to the coprimeness conditions of (14), the integers $M_1,M_2,M_3$, and $M$ are actually mutually relatively prime, thus each of them must be an integer square.  We have,
\setcounter{equation}{15}

\begin{equation}
\left. \begin{array}{l}
M_1 =k^2_1,\ M_2 = k^2_2,\ M_2 = k^2_3,\ M=k^2\\
\\
{\rm for\ positive\ integers}\ k_1,k_2,k_3,k\ {\rm satisfying}\\
\\
(k_1,k_2) = (k_1,k_3) = (k_2,k_3) = 1 = (k,k_1k_2k_3).
\end{array}\right\} \label{E16}
\end{equation}

\noindent Applying the four formulas  in (14) and those in (\ref{E16}), and combining these with the equation ${\displaystyle \frac{N_1}{d}+\frac{N_2}{d} + \frac{N_3}{d} = \frac{N}{d}}$ we arrive at,

\vspace{.15in}
\framebox{ $D_1d_{12}d_{13}k^2_1 + D_2d_{12}d_{23} k^2_2 + D_3 d_{13}d_{23}k^2_3 = D_1D_2D_3 \cdot k^2$ \hspace{.80in} (17)}

\setcounter{equation}{17}

\vspace{.15in}

\noindent We need two Lemmas:

\vspace{.15in}

\noindent{\bf Lemma 1:}  If $A,B,C$ are odd integers, then $AB+BC+AC \equiv3$(mod 4).

\vspace{.15in}

\noindent {\bf Proof:}  Apply the identity $(A+B+C)^2 = A^2+B^2+C^2 +2(AB+AC+BC)$.  Since $A+B+C\equiv 1$(mod 2), we have $(A+B+C)^2 \equiv 1$(mod 8). Hence, from $(A+B+C)^2 \equiv  A^2 \equiv B^2 \equiv C^2 \equiv 1$(mod 8) and the above identity we obtain $(A+B+C)^2 - (A^2+B^2+C^2)\equiv 2(AB+AC+BC)$(mod 8) $\Rightarrow 1-3 \equiv 2(AB+AC+BC)$(mod 8);

\noindent $6 \equiv 2 (AB+AC+BC)$(mod 8) $\Rightarrow AB+BC+AC \equiv 3$(mod 4)

\vspace{.15in}

\noindent{\bf Lemma 2:}  In equation (17), at least one of the following nine integers $D_1,D_2,D_3,d_{12},d_{13},d_{23},k_1,k_2,k_3$, must be even.

\vspace{.15in}

\noindent{\bf Proof:}  If all these nine integers were odd, then so would the integer $k$, as a congruence modulo 2 in (17) shows.  We would have $k^2_1\equiv k^2_2\equiv k^2_3\equiv k^2\equiv1$(mod 4) and if we multiply both sides of (17) by the product $D_1D_2D_3$ and also make use of $D^2_1D^2_2D^2_3 \equiv 1 \equiv D^2_1 \equiv D^2_2 \equiv D^2_3$(mod 4), we end up with the congruence,

\begin{equation}
D_2D_3d_{12}d_{13} + D_1D_3d_{12}d_{23} + D_1D_2 d_{13}d_{23} \equiv {\rm (mod 4)} \label{E18}
\end{equation}

\noindent However, by applying Lemma 1, with $A = D_1D_2d_{12}, B=D_2D_3d_{23}$, and $C=D_1D_3d_{13}$ we obtain

$$
D^2_2D_1D_3d_{12}d_{23} + D^2_1 D_2 D_3 d_{12} d_{13}+ D_1 D_2 d_{13} d_{23} \equiv 3 {\rm (mod\ 4)},
$$

\noindent And since $D^2_1 \equiv D^2_2 \equiv D^2_3 \equiv 1$(mod 4), the last congruence implies

$$
D_1D_3d_{12}d_{23} + D_2D_3d_{12}d_{13} + D_1D_2d_{13}d_{23} \equiv 3\rm{(mod 4)},
$$

\noindent contradicting congruence ({\ref{E18}).

We now go back to the four formulas found in (14), combining them with (\ref{E16}) and (\ref{E3}) to obtain the following set of formulas:

\vspace{.15in}

\framebox{$ \begin{array}{cr}
{\displaystyle \frac{N_1}{d}}= D_1d_{12}d_{13}k^2_1,\ {\displaystyle \frac{N_2}{d}} = D_2d_{12}d_{23}k^2_2,\ {\displaystyle \frac{N_3}{d}} = D_3 d_{13}d_{23}k^2_3,&\\
\\
{\displaystyle \frac{N}{d}} = D_1D_2D_3k^2;\ {\rm and\ area}\ A = \sqrt{{\displaystyle \frac{N\cdot N_1N_2N_3}{16}}} ; &\\
\\
4A = D_1D_2D_3d_{12}d_{23} d_{13}k_1k_2k_3k\cdot d^2;&\\
\\
{\rm and\ sides}\ a,b,c\ {\rm given\ by},& \hspace*{.6in}(19)\\
\\
a = {\displaystyle \frac{d\cdot d_{23}\cdot (D_2d_{12}k^2_2 + D_3d_{13}k^2_3)}{2}}& \\
\\
b = {\displaystyle \frac{d\cdot d_{13}\cdot (D_1d_{12}k^2_1+ D_3d_{23}k^2_3)}{2}}& \\
\\
c = {\displaystyle \frac{d\cdot d_{12} \cdot (D_1d_{13}k^2_1 + D_2d_{23}k^2_2)}{2}} & 
\end{array}$}

\setcounter{equation}{19}
\vspace{.15in}

\noindent Now, apply Lemma 2.  We know that at \underline{least one} of the numbers $k_1,k_2,k_3,d_{12}$,\\
$d_{13},d_{23}, D_1,D_2, D_3$, must be even.  This clearly shows in (19) that one or two of ${\displaystyle \frac{N_1}{d}, \frac{N_2}{d},\frac{N_3}{d}}$ will be even, $\left(\!\!{\rm it\,can\, not\,be\, all\, of\, them\, because}\,\left({\displaystyle \frac{N_1}{d},}\right.\right.$ $\left.{\displaystyle \frac{N_2}{d},\frac{N_3}{d}} \right)$
$\left.\begin{array}{c}\\ = 1\\ \\\end{array} \right)$.  This implies that $d$ must be an even integer.  Why?  Because if $d$ were odd, then one or two of the integers $N_1,N_2,N_3$ would be even; but the formulas in (3) would then imply that not all $a,b,c$ are integers, contrary to the central assumption of this paper, namely that $a,b,c,$ are the integer side lengths of a triangle whose area is rational.

\noindent Hence, \framebox{$d$ must be even}.  Since $d$ is even, and at least one of $k_1,k_2,k_3,d_{12},$\newline
$d_{13},d_{23}, D_1,D_2 D_3$ is even, the formula that gives the area $A$ in (19) implies that the 

\noindent\framebox{area $A$ is in fact an even integer}.

\vspace{.15in}

\noindent \framebox{\begin{tabular}{l} {\bf Note 1:}  If we take a close look at the coprimeness conditions in (14),\\ 
the formulas in (16), and equation (17), we see that at most two of the\\
nine integers
$k_1,k_2,k_3,d_{12},d_{13}, d_{23}, D_1,D_2,D_3$, can be even.  In fact \\
if two of these nine integers are even, then one is $d_{i_1i_2}$ and the other\\
$k_{i_1}$; or alternatively one is $D_{i_1}$ and the other $k_{i_1}$ where $(i_1,i_2,i_3)$ \\
is a permutation of $(1,2,3)$.
\end{tabular}}

\vspace{.15in}

As the listing of all (not exceeding 999) triangle area numbers (at the end of this section) shows, some area numbers are multiples of 4 and some are congruent to 2 modulo 4.  But as it turns out, all of them are divisible by 3.  We have the following.

\vspace{.15in}

\noindent \framebox{\begin{tabular}{l}
{\bf Result 3}\\
(a)  {\it Every triangle area number} $A$ {\it is an even integer}.\\
(b) {\it If} $A \equiv 2$({\it mod 4}), {\it then either exactly one of} $d_{12},d_{13},d_{23}$\\ 
{\it is congruent to 2 modulo 4; the other two are both odd;  all the integers}\\
$k_1,k_2,k_3,k,D_1,D_2,D_3$ {\it are odd; and} $d \equiv 2$({\it mod 4}); \\
{\it or alternatively, exactly one of} $D_1,D_2,D_3$ {\it is congruent to 2}({\it mod 4}); \\
{\it the other two are both odd; and} $k_1,k_2,k_3,k,d_{12},d_{13},d_{23}$\\
{\it are odd; and} $d \equiv 2$({\it mod 4}).\\
(c) {\it Each triangle area number} $A$ {\it is divisible by} $3$.
\end{tabular}}

\vspace{.15in}

\noindent{\bf Proof:}

\begin{enumerate}
\item[(a)] This has already been shown and explained above (see explanation below (19)).
\item[(b)] Suppose $A \equiv2$(mod 4).  Then by the formula for $4A$ in (19), the fact that $d$ is even, and Lemma 2, it follows that $d \equiv 2$(mod 4) and that exactly one of the nine integers $d_{12},d_{13},d_{23},D_1,D_2,D_3,k_1,k_2,k_3$, must be congruent to 2 modulo 4; the other eight of them must all be odd, as well as the integer $k$.  To conclude the proof, it suffice to demonstrate that none of the integers $D_1,D_2,D_3,k_1,k_2,k_3$ can be congruent to 2 mod 4.  Suppose that one of $k_1,k_2,$ or $k_3$ is congruent to 2 mod 4; without loss of generality, say $k_1 \equiv 2$(mod 4), and $d_{12},d_{13},d_{23},D_1,D_2,D_3$, and $k$ are all odd.  We have a contradiction modulo 2, since the left-hand side of (17) will be an even integer, but the right-hand side will be odd.
\item[(c)]  To establish that $A$ must be divisible by $3$, we show that (at least) one of the ten integers $d_{12},d_{13},d_{23},D_1,D_2,D_3,k_1,k_2,k_3,k$ must be divisible by 3; it would then follow by the formula for $4A$ in (19), that $A$ must be a multiple of $3$.  We argue by contradiction; let us suppose that none of the ten numbers $d_{12},d_{13},d_{23},D_1,D_2,D_3,k_1,k_2,k_3,k$,  is divisible by 3.  We have $k^2_1 \equiv k^2_2\equiv k^2_3 \equiv k^2 \equiv 1$(mod 3); applying this to (17) considered modulo 3 we obtain,

\begin{equation}
D_1D_{12}d_{13} + D_2d_{12}d_{23} + D_3 d_{13}d_{23} \equiv D_1D_2D_3{\rm (mod 3)}
\label{E20}
\end{equation}

Since $(d_{12}d_{13} ) \cdot (d_{12}d_{13}) \cdot (d_{12}d_{13}) = d^2_{12}d^2_{13}d^2_{23} \equiv 1$(mod 3);
\begin{enumerate}
\item[(i)] either each of the three products $d_{12}d_{13},d_{12}d_{23},d_{13}d_{23}$, is congruent ot 1 modulo 3 or 
\item[(ii)] two of these three products are congruent to $-1$ modulo 3, the third congruent to 1 modulo 3.
\end{enumerate}
\end{enumerate}

\noindent If (i) is the case then (\ref{E20}) implies, $D_1 + D_2+E_3\equiv D_1D_2D_3$(mod 3).  Each of $D_1,D_2,D_3$ can be either $1$ or $-1$ (mod 3).  No combination of nonzero values modulo 3 of $D_1,D_2,D_3$ can satisfy the latter congruence; it is impossible.

\noindent If (ii) is the case, we can take, without loss of generality, $d_{12}d_{13}\equiv 1$ and $d_{12}d_{23} \equiv d_{13}d_{23} \equiv -1$(mod 3). Using (\ref{E21}) we obtain $D_1-D_2-D_3 \equiv D_1D_2D_3$(mod 3).  Again this has no solution with $D_1D_2D_3 \not\equiv 0$(mod 3).  To see this more clearly consider the fact that at least two of $D_1D_2D_3$ must have the same value mod 3:  $D_1\equiv D_2$ or $D_1D_2D_3$ or $D_2 \equiv D_3$ (the ``or" is not exclusive here); each of these three possibilities combined with $D_1-D_2-D_3 \equiv D_1D_2D_3$(mod 3), implies $D_1D_2D_3 \equiv 0$(mod 3), in violation of $D_1D_2D_3 \not\equiv 0$(mod 3); we omit the details. \hfill End of proof.  $\square$

\noindent We now turn our attention to (17):  If each of the integers $D_1d_{12}d_{13},D_2d_{12}d_{23},$ $D_3d_{13}d_{23}$, and $D_1D_2D_3$ is an integer square, equation (17) procudes a solution to (\ref{E1}); also, according to \underline{Note 2} (above Result 3), the six integers $d_{12},d_{13},d_{23},D_1,D_2,D_3$ are mutually relatively prime; thus $D_1d_{12}d_{13}, D_2d_{12}d_{23},$ $D_3d_{13}d_{23}$, and $D_1D_2D_3$ will all be the integer squares if and only if each integer $d_{12},d_{13},d_{23}, D_1,D_2,D_3$ is an integer square.

\noindent We set

\begin{equation}
D_1 = L^2_1,\ D_2 = L^2_2,\ D_3 = L^2_3,\ d_{12} = \delta ^2_{12}, \ d_{13} = \delta^2_{13},\ d_{23} = \delta ^2_{23}
\label{E21}
\end{equation}

\noindent with $L_1,L_2,L_3,\delta_{12},\delta_{13},\delta_{23}$ mutually relatively prime.  By (\ref{E21}) and (17) we obtain,

\begin{equation}
(L_1\cdot \delta_{12}\cdot \delta_{13}\cdot k_1)^2 + (L_2\cdot \delta_{12}\cdot \delta_{23}\cdot k_2)^2 + (L_3 \cdot \delta_{13}\cdot \delta_{23}\cdot k_3)^2 = (L_1L_2L_3k)^2 
\label{E22}
\end{equation}

\noindent Now it is clear if $(x_0,y_0,z_0,t_0)$ is a solution to (\ref{E1}), we must have (up to symmetry) according to (\ref{E22}),

\begin{equation}
x_0 = L_1\delta_{12}\delta_{13}k_1,\ \ y_0 = L_2\delta_{12}\delta_{23} k_2,\ \ z_0 = L_3 \delta_{13}\delta_{23}k_3,\ \ t_0 = L_1L_2L_3 k.
\label{E23}
\end{equation}

\noindent By considering the coprimeness conditions of (14) and \underline{Note 2} (above Result 3), we can make the following observation:

\begin{enumerate}
\item[(i)] The solution $(x_0,y_0,z_0,t_0)$ to (\ref{E1}) must satisfy $(x_0,y_0,z_0) = 1$
\item[(ii)] $\delta_{12} = (x_0,y_0),\ \delta_{13} = (x_0,z_0),\ \delta_{23} = (y_0,z_0)$
\item[(iii)] $L_1 = \left( {\displaystyle \frac{x_0}{(x_0,y_0)\cdot(x_0,z_0)}},t_0\right) ,\ L_2 = \left( {\displaystyle \frac{y_0}{(y_0,x_0)\cdot(y_0,z_0)}},t_0 \right)$, and
$L_3 = \left( {\displaystyle \frac{z_0}{(z_0,x_0)\cdot (z_0,y_0)}},t_0\right)$
\end{enumerate}

\noindent The three conditions (i), (ii), and (iii) clearly shows that given any solution $\{x_0,y_0,z_0,t_0\}$ to (\ref{E1}), with $(x_0,y_0,z_0) = 1$, the integers $\delta_{12},\delta_{13},\delta_{23},L_1,L_2$, and $L_3$ are uniquely determined up to symmetry or permutations; hence by (\ref{E23}), it follows that the $k_1,k_2,k_3,k$ are likewise uniquely determined. If we make use of the area formula (19), combined with (\ref{E21}) and (\ref{E23}), a computation yields 

\begin{equation}
A = {\displaystyle \frac{x_0\cdot y_0 \cdot z_0 \cdot d^2}{4}} \label{E24}
\end{equation}

\noindent Note that according to \underline{Observation 1} of {\bf Section 3}, in any solution $(x,y,z,t)$ to (\ref{E1}), at least two of $x,y,z$ must be even.  thus, since $d$ is even, (\ref{E24}) clearly shows that the area number $A$  must be \underline{divisible by 4}.

\vspace{.15in}

\noindent {\bf Definition 2:}  Let $(x_0,y_0,z_0,t_0)$ be a positive integer solution to (\ref{E1}) that satisfies $(x_0,y_0,z_0)=1$, and $d$ any even positive integer.  The positive integer $A = {\displaystyle \frac{x_0\cdot y_0 \cdot z_0 \cdot d^2}{4}}$ is called {\bf a solid rectangular number}.  

\vspace{.15in}

\noindent It is clear by now that every solid rectangular number arises directly from equation (\ref{E1}); moreover every solid rectangular number is the area of a \underline{unique} triangle (up to the geometric transformations of reflection and symmetry):  Indeed, by using the formulas for the sides $a,b$, and $c$ in (19) in conjunction with (\ref{E21}) and (\ref{E23})  a computation implies,

\begin{equation}
\left\{ \begin{array}{ll}
a = {\displaystyle \frac{d \cdot (y^2_0 + z^3_0)}{2}},\ b = {\displaystyle \frac{d\cdot (x^2_0 +z^2_0)}{2}},\ c = {\displaystyle \frac{d\cdot (x^2_0 +y^2_0)}{2}}\\
\\
{\rm with} \ (x_0, y_0,z_0) = 1
\end{array}\right\}
\label{E25}
\end{equation}

\noindent Note that the formulas in (\ref{E25}) are equivalent with the formulas of \underline{Result 2} in Section 5.  Indeed, since the formulas of (\ref{E25}) hold under condition $(x_0,y_0,z_0) = 1$, they obviously satisfy the formulas of \underline{Result 2}.  \underline{Conversely}, if we factor out the greatest common divisor $(x_0,y_0,z_0)$ of $x_0,y_0,z_0$ in the formulas of \underline{Result 2}, we immediately see that (\ref{E25}) is satisfied.  Hence these formulas in (\ref{E25}) (or in \underline{Result 2}) generate and describe all the triangles whose areas satisfy (\ref{E24}), i.e., they are solid rectangular numbers.  We end this paper by listing all triangle area numbers not exceeding 999; also listing two important subsets; the Pythagorean numbers and the solid rectangular numbers.  We must note here that there are many, in fact infinitely many, triangle area numbers that correspond to different triangles.  It has been know for quite some time that there exist infinitely many pairs of Pythagorean triangles or triples with the same area number (Pythagorean numbers).  In fact, there exist infinitely many triples of Pythagorean triangles with the same area number, explicitly parametrically described, refer to [2] for more details.  There are exactly 96 triangle area, one digit, two digit, and three digit numbers.

\vspace{.25in}

\noindent{\bf Triangle Area Numbers $\leq$ 999:}  6, 12, 24, 30, 36, 42, 48, 54, 60, 66, 72, 84, 90, 96, 108, 114, 120, 126, 132, 144, 150, 156, 168, 180, 192, 198, 204, 210, 216, 234, 240, 252, 264, 270, 288, 294, 300, 306, 324, 330, 336, 360, 378, 384, 390, 296, 408, 420, 432, 456, 462, 468, 480, 486, 504, 510, 522, 528, 540, 546, 570, 576, 588, 594, 600, 624, 630, 648, 660, 672, 684, 690, 714, 720, 726, 744, 750, 756, 768, 780, 792, 798, 810, 816, 840, 864, 876, 900, 924, 930, 936, 960, 966, 972, 984, 990.

\vspace{.25in}

\noindent {\bf Triangle Area Numbers that are also Pythagorean:}  6, 24, 30, 54, 60, 84, 96, 120, 150, 180, 210, 216, 270, 294, 330, 384, 420, 480 486, 504, 540, 546, 600, 630, 726, 750, 840, 864, 924, 960, 990.

\vspace{.25in}

\noindent {\bf Triangle Area Numbers that are also solid Rectangular:}  12, 972.

\vspace{.25in}

\noindent Note that for example, the Pythagorean number 210 corresponds to different right triangles or triples:  the triangle with sides $a=37, b=35$, and $c=12$; and the triangle with sides $a=20,b=21$, and $c=20$; they both have area number 210.  If we look at the list of 96 area triangle numbers not exceeding 999, we will see (by a quick count) that exactly twenty-one of them are congruent to 2 mod 4 (refer to Result 2 for the special conditions that must be satisfied when $A \equiv 2$(mod 4), in addition to the coprimeness conditions (14)).  So when $A \equiv 2$(mod 4), we have the following list:  6, 42, 66, 90, 114, 126, 198, 234, 306, 390, 462, 510, 522, 570, 594, 690, 714, 798, 810, 966, 990.

Below we list, for each value of $A\equiv 2$(mod 4) (only for the first seven values of $A$), the corresponding values of the integers $D_1,D_2,D_3,d_{12},d_{13},d_{23},$ $k_1,k_2,k_3, k$ and $d$, as well as the side lengths $a,b,c$ of the triangle(s) whose area is the given number $A$.  Recall that $D_1,D_2,D_3,d_{12},d_{13},d_{23},k_1,k_2,k_3,k$, must satisfy the coprimeness conditions (14) and equation (17).

\newpage

\begin{center}
{\bf Table}
\end{center}

\noindent $\begin{array}{||c|c|c|c|c|c|c|c|c|c|c|c|c|c|c||}
\hline\hline
D_1 & D_2 & D_3 & d_{12} & d_{13} & d_{23}& k_1 & k_2 & k_3 & k & d & a & b & c & A\\
\hline
1&2&3&1&1&1&1&1&1&1&2&5&4&3&6=2\cdot 3\\ \hline
1&3&7&1&1&2&1&1&1&1&2&20&15& 7 & 42=2\cdot 3\cdot 7\\ \hline
2&1&11&1&1&1&1&3&1&1&2&20&13&11&66=2\cdot 3\cdot 11\\ \hline
1 & 1 & 3 & 2 & 1 & 5 & 1 & 1& 1& 3 & 2 & 25 & 17 & 12 & 90=2\cdot 3^2\cdot 5\\ \hline
1 & 1&19 & 1&1&1&1&3&1&1&2&37&20&19 & 114=2\cdot 3\cdot 19\\ \hline
2 & 3 & 1 & 1&1&1&1&1&7&3&2&52&51&3&126=2\cdot 3^2 \cdot 7\\ \hline
1&11&6&1&1&1&1&1&3&1&2&65&55&12& 198 = 2\cdot 3^2 \cdot 11 \\ \hline\hline
\end{array}$

\vspace{.15in}

\noindent \underline{\bf  Prove or Disprove}

\vspace{.15in}

\noindent \underline{\bf Conjecture:}  No Pythagorean number is a solid rectangular number.

\begin{center}
\underline{\bf AMS Classification Numbers}
\end{center}

\vspace{.15in}

\noindent 11A99, 11A25, 11D9.


\begin{thebibliography}{99}
\bibitem[Sierpinski] Sierpinski, W., \underline{Elementary  Theory of Numbers}, p. 67, Warsaw, 1964.
\bibitem[Beiler] Bieler, A. H., Recreations in the Theory of Numbers, Dover Publications, Inc., New York, 1966.
\end{thebibliography}
\end{document}